\documentclass[letterpaper, 10pt, conference]{ieeeconf}      

\IEEEoverridecommandlockouts                              
\usepackage{graphics} 
\usepackage{epsfig} 
\usepackage{mathptmx} 
\usepackage{times} 
\usepackage{amsmath} 
\usepackage{amssymb}  

\usepackage[mathcal]{euscript} 
\usepackage{mathtools}
\usepackage{mathrsfs}
\usepackage{dsfont}
\usepackage[noadjust]{cite}

\usepackage{graphicx}
\graphicspath{{Figures/}} 
\usepackage{caption} 
\usepackage{subfig} 

\usepackage{multirow}

\usepackage{algorithm2e}
\SetKwInput{KwInput}{Input}

\usepackage{url}

\usepackage{xcolor}
\definecolor{USred}{rgb}{0.74,0.1,0.1}
\definecolor{USblue}{rgb}{0.2,0.2,0.7}
\definecolor{green1}{cmyk}{0.82,0,1,0.3}

\newcommand{\numberset}{\mathbb}
\newcommand{\NN}{\numberset{N}}

\newcommand{\RR}{\numberset{R}}

\newcommand{\ZZ}{\numberset{Z}}
\newcommand{\EE}{\mathbb{E}}

\newcommand{\A}{\mathcal{A}}
\newcommand{\B}{\mathcal{B}}
\newcommand{\C}{\mathcal{C}}

\newcommand{\R}{\mathcal{R}}
\newcommand{\T}{\mathcal{T}}

\newcommand{\U}{\mathcal{U}}
\newcommand{\V}{\mathcal{V}}

\newcommand{\M}{\mathcal{M}}

\newcommand{\Q}{\mathcal{Q}}

\newcommand{\sO}{\mathcal{O}}
\newcommand{\sH}{\mathcal{H}}

\DeclareMathOperator{\col}{col}
\DeclareMathOperator{\1}{\mathds{1}} 

\DeclareRobustCommand{\vect}[1]{
	\ifcat#1\relax
	\boldsymbol{#1}
	\else
	\mathbf{#1}
	\fi}

\newtheorem{lemma}{Lemma}

\title{An Extended Kalman Filter for Data-enabled Predictive Control}

\author{Daniele Alpago, Florian D\"orfler and John Lygeros 
\thanks{Research supported by the ERC under the project OCAL, the SEMP and the Gini Foundation. D. Alpago is with the Department of Information Engineering at University of Padova, Italy, {\tt\small dalpago@dei.unipd.it}; F. D\"orfler and J. Lygeros are with the Department of Information Technology and Electrical Engineering at ETH Z\"urich, Switzerland, {\tt\small dorfler@ethz.ch}, {\tt\small jlygeros@ethz.ch}.}%
}

\begin{document}

\maketitle
\thispagestyle{empty}
\pagestyle{empty}

\begin{abstract}
The literature dealing with data-driven analysis and control problems has significantly grown in the recent years. Most of the recent literature deals with linear time-invariant systems in which the uncertainty (if any) is assumed to be deterministic and bounded; relatively little attention has been devoted to stochastic linear time-invariant systems. As a first step in this direction, we propose to equip the recently introduced Data-enabled Predictive Control algorithm with a data-based Extended Kalman Filter to make use of additional available input-output data for reducing the effect of noise, without increasing the computational load of the optimization procedure.
\end{abstract}

\section{Introduction}\label{sec:intro}
With the increasing complexity of applications in science and engineering, the use of model-based control techniques is becoming more and more challenging as they usually require accurate descriptions of model and uncertainties, often difficult to obtain. The need of data-driven approaches was already perceived in the 1980s with the advent of system identification \cite{Lju99} and adaptive control \cite{Goo09}, and received new impetuous in the 1990s with approaches such as Iterative Feedback Tuning, Correlation-based Tuning, Virtual Reference Feedback Tuning \cite{Hjaea98,Karea02,Camea02} and many others. In recent years, data-driven analysis and control has experienced a renewed interest \cite{ZhoZhu13}. One promising approach is due to the rediscovery of a result originally formulated in the context of behavioral system theory \cite{Wil91,WilPol13} by J. C. Willems and coauthors, known as the \emph{Fundamental Lemma} \cite{Wilea05}. The result states that if the input signal is sufficiently rich, all possible trajectories of a deterministic, linear time-invariant (LTI) system can be generated from linear combinations of past trajectories of the same system. This allows one to use a Hankel matrix constructed from input-output data as an implicit representation of the underlying dynamics. This representation has been first exploited for data-driven simulation and linear-quadratic tracking for deterministic LTI systems in a behavioral setting \cite{MarRap08,MarRap07}. This has resulted in a growing stream of literature dealing with data-driven analysis \cite{Waaea20,Romea19} and control problems \cite{Berea19b,BerAll19,DepTes19,Waaea19}.

In this context, particular attention has been devoted to the problem of optimal trajectory tracking \cite{Favea94,Couea19a,Couea19b,Huaea19a,Berea19c}, widely studied in model-based control. The aim is to compute an optimal control policy based on output feedback that drives the system along an output reference while minimizing a stage cost and satisfying input and output constraints. Model Predictive Control (MPC) has been one of the most effective methods to tackle such problems. MPC requires an accurate model of the system and an accurate description of the uncertainty (if any), which might be challenging and expensive to obtain in many applications \cite{Hja05,Bulea00}. The success of MPC and the difficulties arising from the requirement for models have led to the introduction of a Data-enabled Predictive Control (DeePC) algorithm. The latter does not rely on a particular parametric system representation, but rather on a Hankel-matrix representation of the underlying dynamics, constructed from the system's input-output trajectories \cite{Couea19a,Couea19b,Huaea19a}. Unlike classical model-based predictive control methods, DeePC performs the control computation without identifying the system by solving a (parametric) convex optimization problem that includes the Hankel matrix in the constraints. The complexity of such optimization depends on the system dimensions, the MPC horizon and the amount of available data. Therefore, it is not obvious how to efficiently include additional data (for example, past measurements of the output of the closed-loop system) in the procedure, without increasing the computational burden of the optimization that has to be solved online.

Unlike earlier methods, recent data-driven approaches based on the Fundamental Lemma have devoted little attention to stochastic systems. Some approaches adopt a ``robust control'' perspective and treat the uncertainty as a deterministic and bounded sequence \cite{Berea19b,Berea19c,DepTes19}, sometimes affecting just the output. The focus of the present paper is to extend the DeePC algorithm to stochastic LTI systems and design methods for including more data to improve closed-loop performance, without increasing the computational load of the optimization. This is accomplished through a combination of an offline averaging of Hankel matrix predictors and an online, data-driven Extended Kalman Filter (EKF). The latter is based on an implicit model constructed from the parametric solution of the DeePC optimization program. The combination of off-line averaging with online EKF dramatically improves the closed-loop performance, as evidenced by numerical experiments.

In Section \ref{sec:prel} we introduce the notation and recall some preliminary results. In Section \ref{sec:deepc} the DeePC algorithm is reviewed and the issues that motivate our approach are highlighted. Section \ref{sec:ekf} introduces the proposed approach and Section \ref{sec:num} presents numerical experiments. Finally, in Section \ref{sec:concl} we draw conclusions and outline future lines of research.

\section{Notation and Preliminaries} \label{sec:prel}
We recall the definition of persistently exciting signals and a state-space version of a result from behavioral system theory \cite{Wil91,WilPol13}, known as the \emph{Fundamental Lemma} \cite{Wilea05,MarRap08}.%

\vspace{2mm}

\emph{Notation.} Given a matrix $A\in\RR^{m\times n}$, $A^\top$ denotes its transpose and $A^{-1}$ its inverse (if $m=n$). The notation $A(i:j,:)$, $i\le j\le m$, stands for the sub-matrix of $A$ that goes from the $i$-th row of $A$ to the $j$-th row of $A$, included. If $i=j$ we simply write $A(i,:)$. We denote with $I_m$ the identity matrix of dimension $m\times m$. With $\text{diag}(A_1,\dots,A_n)$ we denote the block-diagonal matrix formed with the matrices $A_1,\dots,A_n$. Given a sequence of matrices $A_h,\dots,A_{h+k}$ in $\RR^{m\times n}$, with $k>0$, we denote by $A_{h,h+k}:=\col(A_h,\dots,A_{h+k}):=[A_h^\top\,\cdots\,A_{h+k}^\top]^\top\in\RR^{(k+1)m\times n}$. Given a vector $w\in\RR^m$ we introduce the quadratic form $\|w\|_P^2=w^\top P w$ which coincides with the squared Euclidean norm $\|w\|_2^2$ when $P=I_m$. The symbol $\EE[\cdot]$ denotes the expectation operator.

\vspace{2mm}

\emph{Persistency of Excitation.} Let $w_0,\,\dots,w_{N-1}$ be $N\in\NN$ samples of a given sequence $(w_k)_{k\in\ZZ}$ taking values in $\RR^q$. 
For $M\in\NN$, $M\le N$, let
\[
   \sH_M(w_{0,N-1}) :=
   \begin{bmatrix}
   w_0    & w_1    &     & \cdots               &  & w_{N-M}\\
   w_1    & w_2    &     & \cdots               &  & w_{N-M+2}\\
   \vdots & \vdots &     &\reflectbox{$\ddots$} &  & \vdots\\
   w_{M-1} & w_{M} &     & \cdots               &  & w_{N-1}
   \end{bmatrix}
\]
be the block-Hankel matrix associated to the trajectory $w_{0,N-1}$, with $M$ block-rows and the maximal number of columns $N-M+1$. We say that the signal $w_{0,N-1}$ is \emph{persistently exciting} of order $M$ if the Hankel matrix $\sH_M(w_{0,N-1})$ has full row-rank $qM$. This requires the sequence $w_{0,N-1}$ to be rich and long enough, in particular $N\ge(q+1)\,M-1$.%

\vspace{2mm}

\emph{Fundamental Lemma.} Consider a state-space representation of an $n$-dimensional LTI system
\begin{equation}\label{eq:ssm}
	\left\{
	\begin{split}
	x_{k+1} &= A\,x_k + B\,u_k\\
	y_k &= C\,x_k + D\,u_k 
	\end{split}
	\right.
\end{equation}
where $(x_k)_{k\in\ZZ}$ is the $n$-dimensional state-process, $(u_k)_{k\in\ZZ}$ the $m$-dimensional input process, and $(y_k)_{k\in\ZZ}$ is the $p$-dimensional output process. Given an initial condition $x_0\in\RR^n$ and a sequence of inputs $u_{0,k-1}\in\RR^{mk}$, the output of the system can be written as
\begin{equation}\label{eq:output}
	y_{k-1} = CA^{k-1}x_0 + \left[CA^{k-2}B\,\,\cdots\,\,CB\,\,D\right]u_{0,k-1}, \quad k\ge 2,
\end{equation}
and $y_0 = C\,x_0+D\,u_0$. Let $T\in\NN$ and collect $T$-long input-output trajectories $u_{0,T-1},\,y_{0,T-1}$ of the system \eqref{eq:ssm}. For fixed positive integers $N_p,\,N_f\in\NN$ and $k\ge N_p$, we can associate to the vectors
\begin{equation}\label{eq:pastio}
\begin{aligned}
&u_p^{(k)}:=u_{k-N_p,k-1}, \quad y_p^{(k)}:=y_{k-N_p,k-1},\\  
&u_f^{(k)}:=u_{k,k+N_f-1}, \quad y_f^{(k)}:=y_{k,k+N_f-1}, 
\end{aligned}
\end{equation}
the block-Hankel matrices
\begin{align*}
\begin{bmatrix}
U_p \\ U_f
\end{bmatrix}
&= 
\begin{bmatrix}
u_p^{(N_p)} & \cdots & u_p^{(T-N_f-1)}\\
u_f^{(N_p)} & \cdots & u_f^{(T-N_f-1)}
\end{bmatrix}
=:\sH_{N_p+N_f}(u_{0,T-1}),\\
\begin{bmatrix}
Y_p \\ Y_f
\end{bmatrix} &= 
\begin{bmatrix}
y_p^{(N_p)} & \cdots & y_p^{(T-N_f-1)}\\
y_f^{(N_p)} & \cdots & y_f^{(T-N_f-1)}
\end{bmatrix}
=:\sH_{N_p+N_f}(y_{0,T-1}).
\end{align*}
From \eqref{eq:output} we can then construct the $(m+p)(N_p+N_f)\times (T-N_p-N_f+1)$ data matrix
\begin{equation}\label{eq:hankelm}	
\sH :=
\begin{bmatrix}
U_p\\U_f\\Y_p\\Y_f
\end{bmatrix}
=
\begin{bmatrix}
0          & \multirow{2}{*}{$I_{m(N_p+N_f)}$}\\
0          &                                  \\
\sO_p(A,C) & \T_p(B,D)\\
\sO_f(A,C) & \T_f(B,D)\\
\end{bmatrix}
\begin{bmatrix}
X\\U
\end{bmatrix},
\end{equation}
where $U := [u_0\,\,\cdots\,\,u_{T-1}]$ and $X:=[x_0\,\,\cdots\,\,x_{T-1}]$ are the block-Hankel matrices containing the inputs and the corresponding sequence of states, respectively. Here, 
\begin{align*}
\sO_p(A,C) &:= \text{col}(C,\,CA,\,\cdots,\,CA^{N_p-1}),\\
\sO_f(A,C) &:= \text{col}(CA^{N_p},\,CA^{N_p+1},\,\cdots,\,CA^{N_f-1}),
\end{align*}
are observability matrices,
\begin{equation}\label{eq:Toep}
\T_p(B,D) :=
   \begin{bmatrix}
   D 		   & 0 			 & \cdots & 0      & 0      & \cdots & 0\\
   CB 		   & D 	  		 & \cdots & 0      & 0      & \cdots & 0\\
   CAB         & CB 	     & \cdots & 0      & 0      & \vdots & 0\\
   \vdots      & \vdots      & \vdots & \vdots & \vdots & \cdots & 0\\
   CA^{N_p-2}B & CA^{N_p-3}B & \cdots & D      & 0      & \cdots & 0\\
   \end{bmatrix},
\end{equation}
and $\T_f(B,D)$ is defined similarly.
The following fundamental result from \cite{Wilea05} guarantees that the matrix \eqref{eq:hankelm} can be used in place of the parametric representation \eqref{eq:ssm}, as long as the input $u_{0,T-1}$ is persistently exciting.\vspace{2mm}%

\begin{lemma}[Fundamental Lemma]\label{lem:fdl}
	Assume the system \eqref{eq:ssm} to be controllable and the input trajectory $u_{0,T-1}$ to be persistently exciting of order $N_p + N_f + n$. Then, a sequence $\col(u_p,u_f,y_p,y_f)$ is an input-output trajectory of the system \eqref{eq:ssm} if and only if it is in the range space of $\sH$.
\end{lemma}\vspace{2mm}%

Recalling the persistency of excitation condition, $T\ge(m+1)(N_p+N_f+n)+1$ is a necessary condition for $u_{0,T-1}$ to be persistently exciting of order $N_p + N_f + n$. Lemma \ref{lem:fdl} has been originally proven using the behavioral language in \cite[Theorem 1]{Wilea05}. For an equivalent recent state-space proof see \cite{Waaea20} and \cite[Appendix A]{DepTes19}.

\section{Data-Enabled Predictive Control Review}\label{sec:deepc}
We briefly introduce the recently proposed Data-enabled Predictive Control (DeePC) method \cite{Couea19a}, and highlight some related issues when dealing with stochastic LTI systems. This will serve as a motivation for what follows.

Consider the stochastic version of \eqref{eq:ssm}
\begin{equation}\label{eq:ssmstoch}
\left\{
\begin{split}
x_{k+1} &= A\,x_k + B\,u_k + E\,v_k\\
y_k &= C\,x_k + D\,u_k + F\,v_k
\end{split}
\right.
\end{equation}
where $(v_k)_{k\in\ZZ}$ is $p$-dimensional white noise (zero mean and unit-variance). We collect sufficiently long input-output trajectories $u_{0,T-1},\,y_{0,T-1}$, i.e. $T\ge(m+1)(N_p+N_f+n)+1$. The same computations leading to \eqref{eq:hankelm}, lead to the matrix corresponding to the stochastic model \eqref{eq:ssmstoch}:
\begin{equation}\label{eq:hankelmstoch}	
\sH :=
\begin{bmatrix}
U_p\\U_f\\Y_p\\Y_f
\end{bmatrix}
=
\begin{bmatrix}
0          & \multirow{2}{*}{$I_{m(N_p+N_f)}$} & 0 \\
0          &                &  0   \\
\sO_p(A,C) & \T_p(B,D)      & \T_p(E,F)\\
\sO_f(A,C) & \T_f(B,D)      & \T_f(E,F)\\
\end{bmatrix}
\begin{bmatrix}
X\\U\\V
\end{bmatrix},
\end{equation}
where $V := [v_0\,\,\cdots\,\,v_{T-1}]$ is the block-Hankel matrix gathering the noise samples and $\T_p(E,F),\,\T_f(E,F)$ are block-Toeplitz matrices similar to $\T_p(B,D),\,\T_f(B,D)$ in \eqref{eq:Toep}. 

The DeePC algorithm proposed in \cite{Couea19a} attempts to compute an optimal-control action based on past input-output data coming from the unknown system \eqref{eq:ssmstoch}, without previous identification. The control action is selected through an MPC-like optimization problem that allows one to enforce constraints ensuring safety and performance requirements. The previously collected data is directly used on-line in the MPC optimization problem; the predictor is therefore implicit and arises as the outcome of the optimization problem. In particular, at the generic iteration $k\ge N_p$, the DeePC computes optimal control actions by solving
\begin{equation}\label{eq:DeePC}
\begin{aligned}
&\min_{g,u_i,y_i} \quad \sum_{i=1}^{N_f}\|y_i-r_{k+i}\|_Q^2 + \|u_i\|_R^2
                  + \lambda_y\|Y_pg-y_p^{(k)}\|_2^2 + \lambda_g\|g\|_2^2\\
&\text{subject to } \quad  \begin{bmatrix}
U_p\\U_f\\Y_f
\end{bmatrix} g = \begin{bmatrix}u_p^{(k)}\\u_i\\y_i\end{bmatrix},\quad
\begin{split}
u_i\in\U,\quad i=1,\dots, N_f,\\
y_i\in\mathcal{Y},\quad i=1,\dots, N_f,
\end{split}
\end{aligned}
\end{equation}  
where $N_f$ is the prediction horizon, $r\in\RR^{pN_f}$ is the output reference signal we want to track, $\U\subseteq\RR^m,\,\mathcal{Y}\subseteq\RR^p$ are the input and output constraint sets, respectively, $Q\in\RR^{p\times p}$ is the output cost matrix (positive semidefinite), $R\in\RR^{m\times m}$ is the control cost matrix (positive definite), $\lambda_y\ge 0$ and $\lambda_g\ge0$ are the regularization parameters, and $u_p^{(k)},\,y_p^{(k)}$ are the most recent $N_p$ input-output measurements from \eqref{eq:ssmstoch}, according to the notation \eqref{eq:pastio}. For simplicity we consider input and output box constraints of the form $\U=[u_\text{min},\,u_\text{max}]$ and $\mathcal{Y}=[y_\text{min},\,y_\text{max}]$, respectively. Note that the block-Hankel matrices $U_p,U_f,Y_p,Y_f$ are fixed throughout the online iterations.

If we let $g_k^\star$ be the optimal solution of problem \eqref{eq:DeePC} at iteration $k$, DeePC provides an implicit predictor $y_f^{(k)}=Y_fg_k^\star$ whose model is never derived explicitly but whose predictions are implicitly used for obtaining the optimal control actions $u_f^{(k)}=U_fg^\star_k$. Problem \eqref{eq:DeePC} is solved in a receding-horizon fashion: of the $N_f$-long optimal control sequence $u_f^{(k)}=\col(u_k^\star,\dots,u_{k+N_f-1}^\star)$ solution of \eqref{eq:DeePC}, we apply a sub-sequence $u_k^\star,\dots,u_{k+N_\text{c}}^\star$, for some $N_\text{c}\le N_f-1$, to the system, update $u_p^{(k)},\,y_p^{(k)}$ to the most recent input-output measurements and set $k$ to $k+N_\text{c}+1$.

Equation \eqref{eq:DeePC} is a relaxation of the corresponding problem enforcing the constraint $Y_p\,g=y_p^{(k)}$, which was proven to be equivalent to the classical receding-horizon MPC in the case of deterministic LTI systems \cite{Couea19a}. Here, to cope with potential infeasibility due to the disturbances, the constraint is substituted with the least-squares regularization term $\|Y_pg-y_p^{(k)}\|_2^2$. The two-norm regularization on $g$ has been introduced to avoid overfitting and it has been shown to relate to distributional robustness of the method with respect to a range of uncertainties \cite{Couea19b}. One can also see that, besides infeasibility issues, considering stochastic models makes both the implicit predictor model $y_f^{(k)}=Y_fg^\star_k$ and the control actions $u_f^{(k)}=U_fg^\star_k$  defined by the optimization \eqref{eq:DeePC} to be non-linear in the past-data $(y_p^{(k)},\,u_p^{(k)})$ (Section \ref{sec:meth}). 

We conclude this section with a note on the effect of acquiring more data. Classical LTI system identification methods use historical data off-line to compute estimates of the matrices in the system dynamics. The resulting matrices are then used on-line to generate state estimates and perform the prediction \cite{Lju99}. Other prediction methods estimate the linear relation from $u_p^{(k)},u_f^{(k)},\,y_p^{(k)}$ to $y_f^{(k)}$ off-line, and then use it online to carry on the predictions, \cite{OveDem12} and references therein. This kind of prediction architectures naturally possesses a data-compression mechanism. For the classical system identification methods only storage of the estimates of the system matrices is required. This depends on input and output dimensions $m$ and $p$ and on the (guessed) state dimension, but not on the amount $T$ of data nor on the prediction horizon $N_f$. Likewise, the storage requirements for linear predictors depends on input and output dimensions $m$ and $p$, and the horizons $N_p$, $N_f$, but not on the amount $T$ of data. If more data become available one can then use it to improve predictions with no need to store it or use it online. The DeePC approach, on the other hand, requires one to carry all the data at every iteration and the size of the optimization problem \eqref{eq:DeePC} increases both with the amount $T$ of available data and with the parameters $m$, $p$, $N_p$, and $N_f$. Though additional data would in principle also be beneficial for a DeePC controller, it is not clear how one can incorporate it without increasing the on-line computational burden.

\section{Method Description}\label{sec:meth}
We introduce a possible way to effectively incorporate more data in the DeePC framework to reduce the effect of noise in the solution of problem \eqref{eq:DeePC}. We propose to use additional data that may be available off-line to de-noise the block-Hankel matrices $Y_p$ and $Y_f$ and to equip the DeePC with an EKF based on data to handle the noise in the on-line measurements $y_p^{(k)}$.

\subsection{Averaging Data Matrices}\label{sec:avehank}
The output matrices $Y_p$ and $Y_f$ are constructed off-line from the trajectory $y_{0,T-1}$. Let $(u_{0,T-1}^{(1)},y_{0,T-1}^{(1)}),\dots,(u_{0,T-1}^{(N)}y_{0,T-1}^{(N)})$ be available $T$-long input-output trajectories and denote with $x_0^{(i)},\,i=1,\dots N$, the corresponding initial conditions. 
Using those additional data on-line to improve the prediction will lead to an intractable optimization problem. However, we can make use of additional trajectories off-line to construct $N$ different data matrices $\sH^{(1)},\dots,\sH^{(N)}$ defined analogously to \eqref{eq:hankelmstoch}, and average those matrices to obtain
\begin{equation}\label{eq:hankave}
\bar{\sH}_N:=\frac{1}{N}\sum_{i=1}^N\sH^{(i)}
=
\begin{bmatrix}
0          & \multirow{2}{*}{$I_{m(N_p+N_f)}$} & 0 \\
0          &                & 0 \\
\sO_p(A,C) & \T_p(B,D)      & \T_p(E,F)\\
\sO_f(A,C) & \T_f(B,D)      & \T_f(E,F)\\
\end{bmatrix}
\begin{bmatrix}
\bar{X}_N\\\bar{U}_N\\\bar{V}_N
\end{bmatrix}
\end{equation} 
where $\bar{V}_N:=\frac{1}{N}\,\sum_{i=1}^N\,V^{(i)}$
is the average of the matrices $V^{(i)}=[v_0^{(i)}\,\,\cdots\,\,v_{T-1}^{(i)}]$, $i=1,\dots,N$, constructed from the different noise realizations $v_{0,T-1}^{(i)}$ affecting the trajectories $y_{0,T-1}^{(i)}$, and similarly $\bar{X}_N$ and $\bar{U}_N$ contain the averaged state trajectories corresponding to the different initial conditions $x_0^{(i)}$ and input trajectories $u_{0,T-1}^{(i)}$, respectively. Since $(v_k)_{k\in\ZZ}$ is assumed to be white-noise, the Law of Large Numbers guarantees that $\bar{V}_N\to 0$ as $N\to\infty$ \cite{Bil95}. The averaging procedure makes use of additional data to mitigate the effect of noise in the data-driven model, hence reducing the risk of overfitting that would be present if the data was used directly in DeePC. Accordingly, the more matrices are involved in the average, the smaller the value of the regularization parameter $\lambda_g$ that gives the optimal closed-loop cost (see Figure \ref{fig:avecost}, Section \ref{sec:num}).

We recognize that such a method heavily exploits the underlying linear structure of the problem. However, this is meant to be a first attempt to exploit additional data for improving the performance of the algorithm when dealing with stochastic systems, without increasing the dimension of the optimization problem \eqref{eq:mpEKFDeePC} to be solved online. Indeed, thanks to the linear structure and superposition, the sub-matrices $\bar{X}_N$ and $\bar{U}_N$ in the averaged data matrix \eqref{eq:hankave} still represent valid system trajectories corresponding to the average of the initial states $x_0=\frac{1}{N}\sum_ix_0^{(i)}$ and the average of the input sequences $u=\frac{1}{N}\sum_iu^{(i)}$, used in the experiments that generated the data for each $\sH^{(i)}$. The only point that requires attention is ensuring that the resulting average input sequence respects the persistence of excitation requirement. The simplest way to ensure this is to assume that the same persistently exciting input sequence is applied in all cases.

\subsection{An EKF for DeePC}\label{sec:ekf}
The averaging procedure represents a simple way to make use of multiple $T$-long trajectories to denoise off-line the data matrix representing the dynamics in the constraints of \eqref{eq:DeePC}. Numerical evidence suggests that additional denoising of the on-line data $y_p^{(k)}$ that enters the cost of \eqref{eq:DeePC} can lead to a further, significant improvement in performance (Section \ref{sec:num}). In a model based setting, such on-line denoising could be performed by a Kalman filter. In our data-driven setting, however, the classical Kalman filter cannot be applied as it requires a model of the system dynamics. Here we show how the non-linear one-step predictor implicit in \eqref{eq:DeePC} can be used to derive an EKF for integrating past measurements into an implicit ``state estimate"; this can in turn be used to improve the asymptotic performance of the algorithm. 

We introduce a fictitious state vector $z_k := \text{col}(y_{k-Np+1},\,\cdots,\, y_k)\in\RR^{pN_p}$
and keep track of a filtered state-estimate $\hat{z}_{k|k}$ and the corresponding error covariance $P_{k|k}$. At the generic iteration $k\ge N_p$, the EKF-DeePC algorithm uses the estimate $\hat{z}_{k|k}$ by solving the following optimization problem
\begin{equation}\label{eq:EKFDeePC}
\begin{aligned}
\min_{g\in\RR^d} &\quad \sum_{i=1}^{N_f}\|(Y_fg)_i-r_{k+i}\|_Q^2 + \|(U_f g)_i\|_R^2\\ 
&\quad + \lambda_y\|Y_pg-\hat{z}_{k|k}\|_2^2 + \lambda_g\|g\|_2^2\\
\text{subject to } &\quad U_pg=u_p^{(k)},\\
&\quad u_\text{min}\le (U_f g)_i \le u_\text{max},\quad i=1,\dots,N_f,\\
&\quad y_\text{min}\le (Y_f g)_i \le y_\text{max},\quad i=1,\dots,N_f.
\end{aligned}
\end{equation}
The formulation \eqref{eq:EKFDeePC} is obtained from \eqref{eq:DeePC} by substituting the constraints $u_i=U_fg$ and $y_i=Y_fg$ in the cost, leaving $g$ as the only decision variable. The crucial difference with respect to \eqref{eq:DeePC} is that the past data $Y_pg$ used in the implicit predictor is now required to fit the state estimate $\hat{z}_{k|k}$ instead of the $N_p$ most-recent measurements $y_p^{(k)}$. Rewriting \eqref{eq:EKFDeePC} as a multi-parametric quadratic program (mp-QP) in the parameter $\theta_k:=\col(\hat{z}_{k|k},\,u_p^{(k)})$ leads to an explicit relation between of the optimum $g^\star_k$ and the parameter $\theta_k$; we exploit this relation to build the EKF. The mp-QP form of \eqref{eq:EKFDeePC} is
\begin{equation}\label{eq:mpEKFDeePC}
\begin{aligned}
\min_{g\in\RR^d} &\quad \frac{1}{2}\,g^\top P g + (G\theta_k+q_k)^\top g + \theta_k^\top H \theta_k + \frac{1}{2}\vect{r}_k^\top\vect{Q}\vect{r}_k\\
\text{subject to } &\quad U_p\,g = B_\text{eq}\,\theta_k,\quad A_\text{in}\,g\le b_\text{in},
\end{aligned}
\end{equation}
where the inequality constraint defined by 
\begin{align*}
b_\text{in}&:=\col(\1_{mN_f} u_\text{max},-\1_{mN_f} u_\text{min},\1_{pN_f} y_\text{max},-\1_{pN_f} y_\text{min}),\\
A_\text{in} &:= \col(U_f,-U_f,Y_f,-Y_f),
\end{align*}
has to be understood \emph{component-wise}, i.e. $(A_\text{in}\,g)_i \le s_i$ for $i=1,\dots,2(m+p)N_f$. The cost is then defined by the reference signal $\vect{r}_k:=\col(r_{k+1},\cdots,r_{k+N_f})$ and the matrices
\begin{align*}
&P:= Y_f^\top\vect{Q}Y_f+U_f^\top\vect{R}U_f+\lambda_yY_p^\top Y_p+\lambda_g\,I_d,\quad 
q_k := -Y_f^\top\vect{Q}\vect{r}_k\\[1mm]
&\begin{aligned}
&G :=\left[\begin{array}{cc}
-\lambda_y Y_p^\top & 0\end{array}
\right]\\
&B_\text{eq} := \left[\begin{array}{cc}
0 & I_{mN_p}
\end{array}\right]
\end{aligned},\quad
H :=\begin{bmatrix}
(\lambda_y/2)\cdot I_{pN_p} & 0\\
0 & 0
\end{bmatrix}.
\end{align*} 
Here $\1_N:=\col(1,1,\cdots,1)\in\RR^N$, $\vect{Q}=\text{diag}(Q,\dots,Q)$ and $\vect{R}=\text{diag}(R,\dots,R)$. Assuming $P=P^\top>0$ and the KKT-matrix for problem \eqref{eq:mpEKFDeePC} to be positive semi-definite (which is always the case for an MPC problem with input weighting matrix $R>0$ \cite{Bemea02}), the optimizer $g_k^\star$ is a piecewise affine function of the parameters, and can be written as $g_k^\star = \tilde{\A}_k\,\hat{z}_{k|k} + \tilde{\B}_k\,u_p^{(k)} + \tilde{h}_k$.
In particular, the equality constraints in \eqref{eq:mpEKFDeePC} and the noise in the model are responsible of the affine structure while the inequality constraints implies this affine relation to hold just locally, i.e. in a neighborhood of the parameter $\theta_k$, known as the \emph{critical region}. The coefficients $\tilde{\A}_k$, $\tilde{\B}_k$ and $\tilde{h}_k$ coming from the KKT conditions for problem \eqref{eq:mpEKFDeePC}, are therefore region-dependent themselves \cite{Bemea02}: we consider the affine expansion of $g_k^\star$ pertaining to a specific value for the parameter $\theta_k$, fixed by the previous iteration.

The fact that, under suitable assumptions, the piecewise affine relation between the optimizer $g_k^\star$ of \eqref{eq:mpEKFDeePC} and the parameter $\theta_k$ can be derived from the KKT conditions for problem \eqref{eq:mpEKFDeePC}, might suggest that the predictor implicit in \eqref{eq:EKFDeePC} could be made explicit and be constructed off-line. Doing this would require one to construct all the regions on which the affine expansion is defined, for all the possible values of the parameters. As the number of these regions scales exponentially with the QP size (parameters plus constraints), hence with the amount $T$ of available data and the horizons $N_p$ and $N_f$, this approach is likely to be computationally intractable. This is the main reason for keeping the predictor implicit.

\begin{figure}[h!]
	\centering
	\includegraphics[scale=0.6]{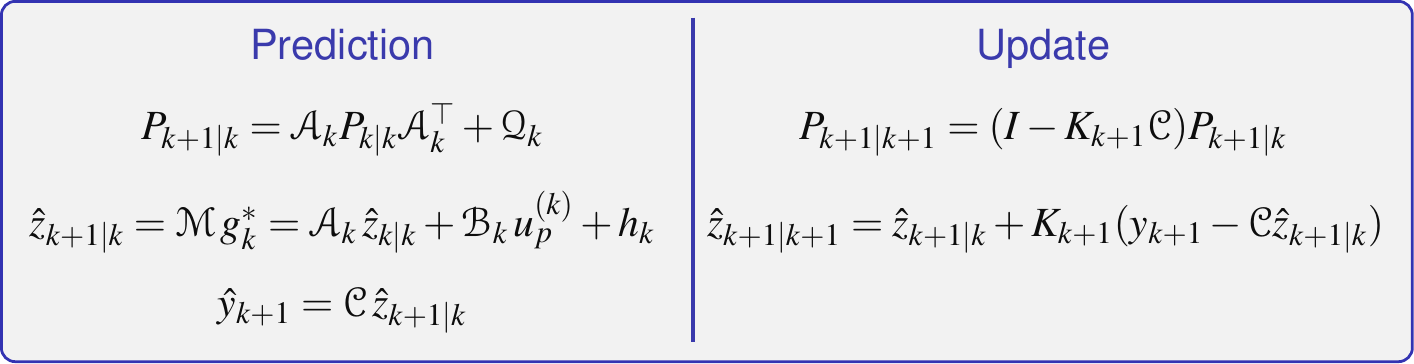}
	\caption{EKF equations. Here, $\C:=[0\,\,\cdots\,\,0\,\,I_p]$ by construction, $\A_k:=\M\,\tilde{\A}_k$, $\B_k := \M\,\tilde{B}_k$, $h_k := \M\,\tilde{h}_k$ and $S_{k+1} = \C P_{k+1|k}\C^\top + \R_k$, $K_{k+1}=P_{k+1|k}\C^\top S_{k+1}^{-1}$, are the variance of the innovation process and the filter gain, respectively.}
	\label{eq:ekf}
\end{figure}

Let $\M:=\text{col}(Y_p(2:pN_p,:),\,Y_f(1,:))$ be the prediction map.
We can exploit the piece-wise affine form of the optimizer $g_k^\star$ to incorporate the implicit predictor provided by \eqref{eq:mpEKFDeePC} in an EKF-like architecture. To fix ideas, suppose we start running the algorithm at $k=N_p$. From an initial guess of the mean $\hat{z}_{N_p|N_p}=\EE[z_{N_p}]$ and the covariance matrix $P_{N_p|N_p}=\EE[(z_{N_p}-\hat{z}_{N_p|N_p})(z_{N_p}-\hat{z}_{N_p|N_p})^\top]$, we then
compute recursively (at every point in time) the standard Kalman filter update steps in Figure \ref{eq:ekf}. We note that, because of its dependence on $g_k^\star$, the implicit predictor is piecewise affine (the matrices $\A_k,\,\B_k,\,h_k$ at the current iteration, depend on the critical region) making this a data-driven analogue to an {\em Extended} Kalman Filter.
\begin{figure}[h!]
	\centering
	\includegraphics[scale=0.7]{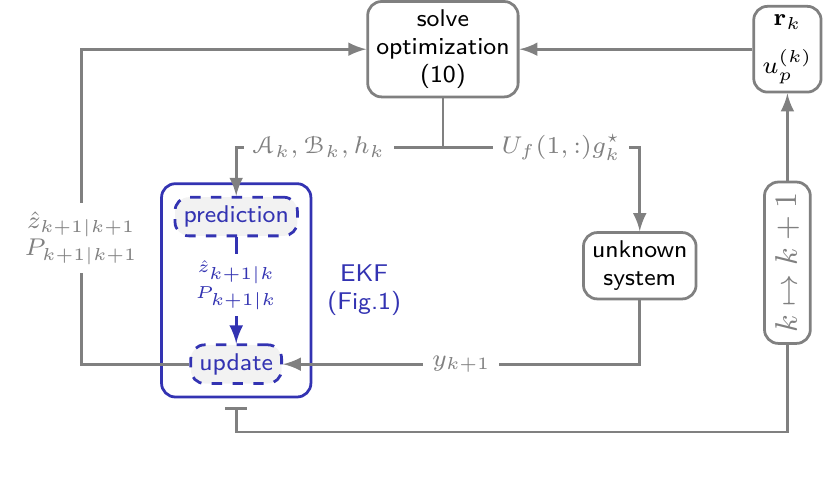}
	\caption{EKF-DeePC Algorithm.}
	\label{fig:ekfdpcalg}
\end{figure}
\vspace{-2mm}

The matrices $\Q_k\in\RR^{pN_p\times pN_p}$ and $\R_k\in\RR^{p\times p}$ are the variances of the process and measurement noise implicitly affecting the state and output dynamics $z_{k+1} = \A_k\,z_k + \B_k\,u_p^{(k)} + h_k$ and $y_k = \C\,z_{k}$, respectively.
In practice, to implement the filter in Figure \ref{eq:ekf} one needs to address the same issues as for any Kalman filter, i.e. choosing the initial conditions $\hat{z}_{N_p|N_p}$ and $P_{N_p|N_p}$, and obtaining an estimate (or guess) of the variances $\Q_k$ and $\R_k$. For a model-based Kalman filter, one can use the residuals of the system identification to obtain estimates for these quantities; we speculate that something similar would be possible using the residuals between the average and individual data matrices in Section \ref{sec:avehank}, though we do not pursue this direction here, due to space limitations. The steps outlined above are summarized in Figure \ref{fig:ekfdpcalg}.

\section{Numerical Validation}\label{sec:num}
We present numerical experiments assessing the effects of the proposed solution in tackling the noise acting in the optimization \eqref{eq:DeePC}. Consider the stochastic system 
\begin{equation}\label{eq:ssmsim}
\left\{
\begin{aligned}
x_{k+1} &= Ax_k + Bu_k + w_k\\
y_k &= Cx_k + v_k
\end{aligned}
\right.,
\quad
A=\begin{bmatrix}
0.8 & 1\\
0 & 0.8
\end{bmatrix},\quad
\begin{aligned}
B &= [0\,\,1]^\top,\\
C &= [1\,\,1].
\end{aligned}
\end{equation}
where $(w_k)_{k\in\ZZ},\,(v_k)_{k\in\ZZ}$ are independent zero-mean Gaussian white noises with covariance $\Sigma_w:=EE^\top$ and $\Sigma_v:=FF^\top$ respectively; below we consider $\Sigma_w=\sigma_w^2 I_n$ and $\Sigma_v = \sigma_v^2 I_p$ and report results for different values of $\sigma_w^2$ and $\sigma_v^2$. Model \eqref{eq:ssmsim} is reachable and observable. The performance metric we will consider throughout is the closed-loop cost $J(u,y) = \sum_{k=1}^{N_\text{sim}}\|y_k-r_k\|_Q^2 + \|u_k\|_R^2,$
where $r_k=5\sin(0.3k)$ is the reference signal. To isolate the effects of averaging and the EKF, for each numerical experiment we tuned the regularization parameters $\lambda_y$ and $\lambda_g$ through exhaustive search to minimize the closed-loop cost for the standard DeePC, the averaged DeePC and the averaged DeePC with EKF. All numerical evidence comes from $100$ repetitions for different data-sets. Unless otherwise stated, we set $T=100$, $N_\text{sim}=100$, $N_p=3$, $N_f=5$, $Q=I_p$, $R=I_m$, $\sigma_w^2=\sigma_v^2=0.5$ and $N=40$ data matrices in the average \eqref{eq:hankave}.

\begin{figure}[h!]
	\centering
	\includegraphics[scale=0.5]{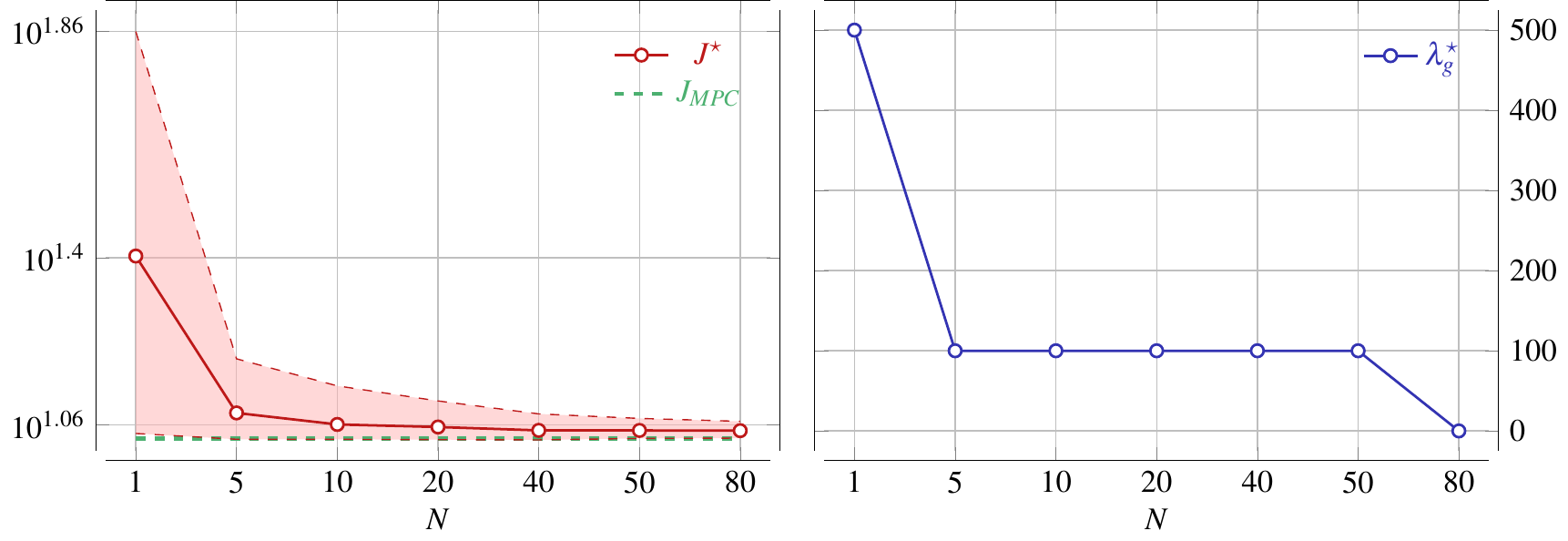}
	\caption{Closed-loop cost (left) and optimal $\lambda_g$ (right) when averaging $N$ data matrices. The dashed-green line shows the MPC cost with the same horizon, but perfect model and full state measurement.}
	\label{fig:avecost}
\end{figure}

\begin{figure}[h!]
	\centering
	\includegraphics[scale=0.5]{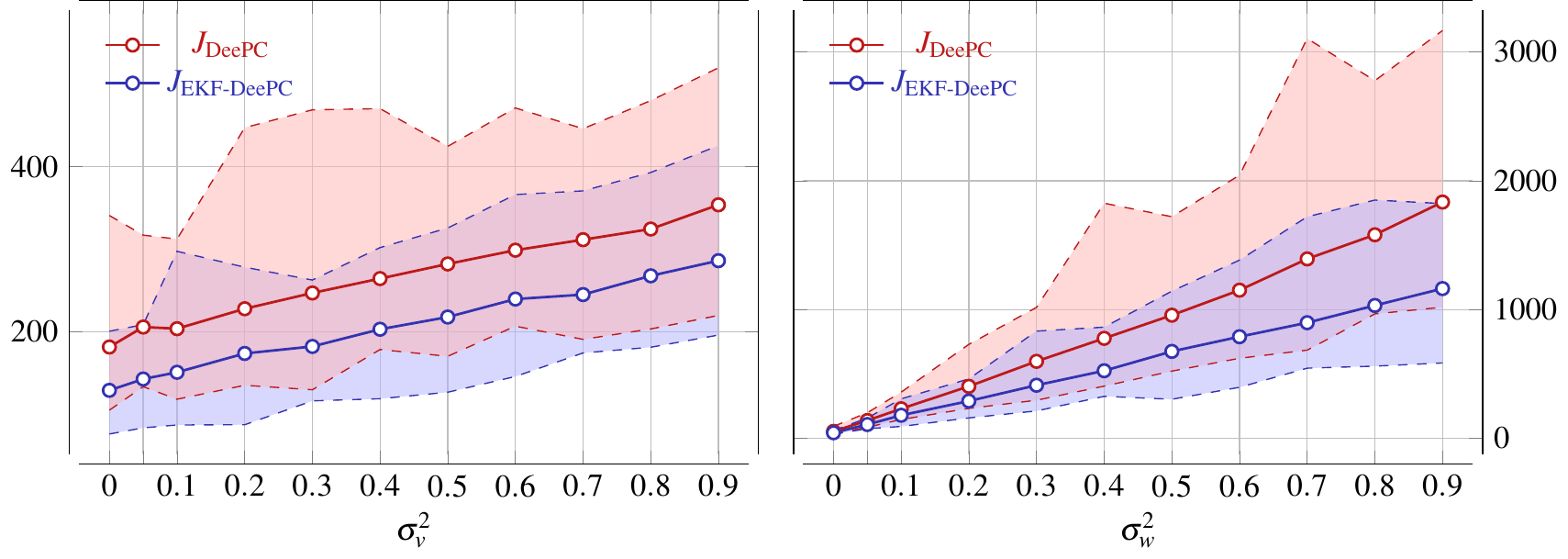}
	\caption{Closed-loop cost of standard DeePC and averaged DeePC with EKF for different $\sigma_v^2$ while keeping $\sigma_w^2=0.1$ (left), and different $\sigma_w^2$ while keeping $\sigma_v^2=0.2$ (right).}
	\label{fig:vars}
\end{figure}

First we show the effect of the averaging introduced in Section \ref{sec:avehank}. To isolate the averaging effect, we consider noisy $Y_p,\,Y_f$ but $y_p^{(k)}$ generated from the model with neither process noise, not measurement noise. The result, shown in Figure \ref{fig:avecost}, is as expected: the cost decreases towards the ideal cost of MPC with the same horizon but perfect model and full state measurement while on the right, the optimal (numerically found) value of the regularization parameter $\lambda_g$ decreases to zero, as we expect from standard results in optimization.
\vspace{-2mm}
Figure \ref{fig:vars} shows a comparison between the closed-loop costs of the standard DeePC and the averaged DeePC with EKF algorithms when varying the noise variances $\sigma_w^2$ and $\sigma_v^2$. Figure \ref{fig:vars} shows that introducing averaging and the EKF substantially improves the performance of the DeePC algorithm. Further experiments (data not shown) with averaged DeePC without EKF and EKF without averaging confirmed that DeePC with a combination of averaging and EKF substantially outperformed all other combinations; indeed this was the case across a broad range of values of the regularizes $\lambda_y$ and $\lambda_g$, suggesting that averaged DeePC with EKF is easier to tune.
\begin{figure}[h!]
	\centering
	\includegraphics[scale=0.5]{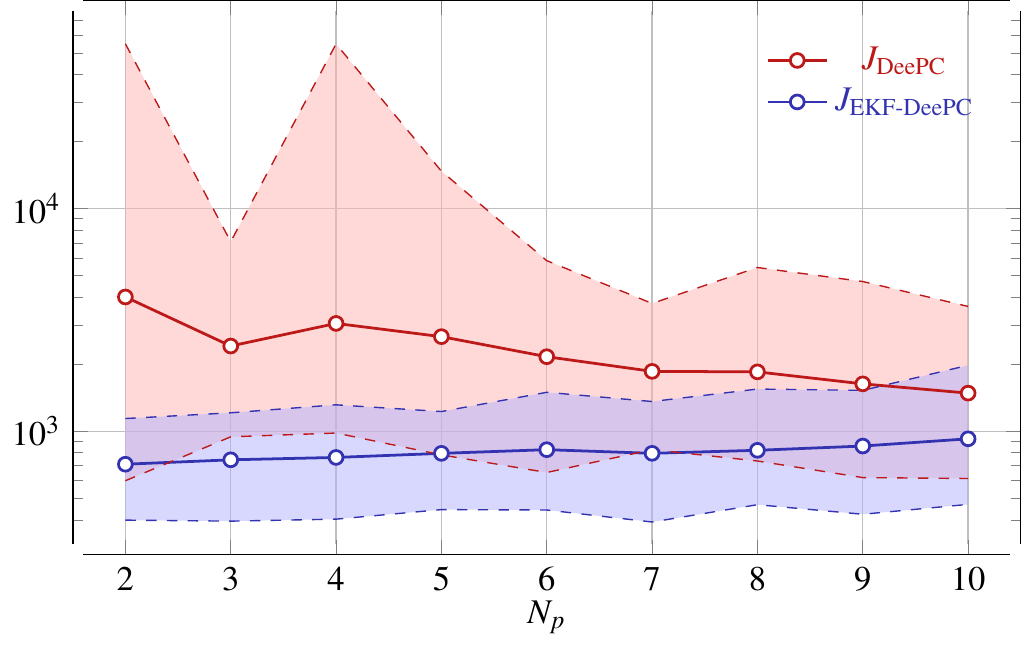}
	\caption{Closed-loop costs while varying the past horizon $N_p$.}
	\label{fig:Npcomp}
\end{figure}
Figure \ref{fig:Npcomp} shows how the closed-loop cost of the two algorithms behaves with respect to the past horizon $N_p$. The improvement in the performance is due to the introduction of the averaging plus the EKF, highlighting the difference between fixed horizon and recursive filtering.

\section{Conclusion and Outlook}\label{sec:concl}
We presented an extension of the data-enabled predictive control (DeePC) algorithm introduced in \cite{Couea19a} to tackle control problems on unknown stochastic LTI systems, by making use of additional data without increasing the dimension of the on-line optimization problem. The procedure features an on-line data-driven EKF that filters out the noise in the measurements, and an (off-line) averaging of multiple data matrices to get a cleaner data-driven model. The performance of the proposed method has been validated experimentally on an LTI stochastic system showing improvements with respect to standard DeePC. Future work includes testing the proposed method on a real-case scenario and comparing it to algorithms such as \cite{Hjaea98,Karea02,Camea02}.

\section*{Acknowledgments}
The authors are grateful to Colin Jones and Francesco Micheli for detailed discussions and the anonymous reviewers for insightful suggestions.

\bibliography{biblio}

\begin{thebibliography}{10}

\bibitem{Lju99}
L.~Ljung, {\em System Identification: Theory for the User}.
\newblock Prentice Hall information and system sciences series, Prentice Hall
  PTR, 1999.

\bibitem{Goo09}
G.~C. Goodwin and K.~S. Sin, {\em Adaptive Filtering Prediction and Control}.
\newblock USA: Dover Publications, Inc., 2009.

\bibitem{Hjaea98}
H.~{Hjalmarsson}, M.~{Gevers}, S.~{Gunnarsson}, and O.~{Lequin}, ``Iterative
  feedback tuning: theory and applications,'' {\em IEEE Control Systems
  Magazine}, vol.~18, no.~4, pp.~26--41, 1998.

\bibitem{Karea02}
A.~Karimi, L.~Mišković, and D.~Bonvin, ``Convergence analysis of an iterative
  correlation-based controller tuning method,'' {\em IFAC Proceedings Volumes},
  vol.~35, no.~1, pp.~413 -- 418, 2002.
\newblock 15th IFAC World Congress.

\bibitem{Camea02}
M.~Campi, A.~Lecchini, and S.~Savaresi, ``Virtual reference feedback tuning: a
  direct method for the design of feedback controllers,'' {\em Automatica},
  vol.~38, no.~8, pp.~1337 -- 1346, 2002.

\bibitem{ZhoZhu13}
Z.-S. Hou and Z.~Wang, ``From model-based control to data-driven control:
  Survey, classification and perspective,'' {\em Information Sciences},
  vol.~235, pp.~3 -- 35, 2013.
\newblock Data-based Control, Decision, Scheduling and Fault Diagnostics.

\bibitem{Wil91}
J.~C. {Willems}, ``Paradigms and puzzles in the theory of dynamical systems,''
  {\em IEEE Transactions on Automatic Control}, vol.~36, pp.~259--294, March
  1991.

\bibitem{WilPol13}
J.~C. Willems and J.~W. Polderman, {\em Introduction to mathematical systems
  theory: a behavioral approach}, vol.~26.
\newblock Springer Science \& Business Media, 2013.

\bibitem{Wilea05}
J.~C. Willems, P.~Rapisarda, I.~Markovsky, and B.~L.~D. Moor, ``A note on
  persistency of excitation,'' {\em Systems \& Control Letters}, vol.~54,
  no.~4, pp.~325 -- 329, 2005.

\bibitem{MarRap08}
I.~Markovsky and P.~Rapisarda, ``Data-driven simulation and control,'' {\em
  International Journal of Control}, vol.~81, no.~12, pp.~1946--1959, 2008.

\bibitem{MarRap07}
I.~Markovsky and P.~Rapisarda, ``On the linear quadratic data-driven control,''
  in {\em 2007 European Control Conference (ECC)}, pp.~5313--5318, IEEE, 2007.

\bibitem{Waaea20}
H.~J. van Waarde, C.~D. Persis, M.~K. Camlibel, and P.~Tesi, ``Willems'
  fundamental lemma for state-space systems and its extension to multiple
  datasets,'' 2020.

\bibitem{Romea19}
A.~{Romer}, J.~{Berberich}, J.~{Köhler}, and F.~{Allgöwer}, ``One-shot
  verification of dissipativity properties from input–output data,'' {\em
  IEEE Control Systems Letters}, vol.~3, pp.~709--714, July 2019.

\bibitem{Berea19b}
J.~Berberich, A.~Romer, C.~W. Scherer, and F.~Allg{\"o}wer, ``Robust
  data-driven state-feedback design,'' {\em arXiv preprint arXiv:1909.04314},
  2019.

\bibitem{BerAll19}
J.~Berberich and F.~Allg{\"o}wer, ``A trajectory-based framework for
  data-driven system analysis and control,'' {\em arXiv preprint
  arXiv:1903.10723}, 2019.

\bibitem{DepTes19}
C.~{De Persis} and P.~{Tesi}, ``Formulas for data-driven control:
  Stabilization, optimality, and robustness,'' {\em IEEE Transactions on
  Automatic Control}, vol.~65, pp.~909--924, March 2020.

\bibitem{Waaea19}
H.~J. van Waarde, J.~Eising, H.~L. Trentelman, and M.~K. Camlibel, ``Data
  informativity: a new perspective on data-driven analysis and control,'' 2019.

\bibitem{Favea94}
W.~Favoreel, B.~D. Moor, and M.~Gevers, ``Spc: Subspace predictive control,''
  {\em IFAC Proceedings Volumes}, vol.~32, no.~2, pp.~4004 -- 4009, 1999.
\newblock 14th IFAC World Congress 1999, Beijing, Chia, 5-9 July.

\bibitem{Couea19a}
J.~Coulson, J.~Lygeros, and F.~D{\"o}rfler, ``Data-enabled predictive control:
  In the shallows of the deepc,'' in {\em 2019 18th European Control Conference
  (ECC)}, pp.~307--312, IEEE, 2019.

\bibitem{Couea19b}
J.~Coulson, J.~Lygeros, and F.~D{\"o}rfler, ``Regularized and distributionally
  robust data-enabled predictive control,'' {\em arXiv preprint
  arXiv:1903.06804}, 2019.

\bibitem{Huaea19a}
L.~Huang, J.~Coulson, J.~Lygeros, and F.~Dorfler, ``Data-enabled predictive
  control for grid-connected power converters,'' {\em arXiv preprint
  arXiv:1903.07339}, 2019.

\bibitem{Berea19c}
J.~Berberich, J.~K{\"o}hler, M.~A. M{\"u}ller, and F.~Allg{\"o}wer,
  ``Data-driven model predictive control with stability and robustness
  guarantees,'' {\em arXiv preprint arXiv:1906.04679}, 2019.

\bibitem{Hja05}
H.~Hjalmarsson, ``From experiment design to closed-loop control,'' {\em
  Automatica}, vol.~41, no.~3, pp.~393 -- 438, 2005.
\newblock Data-Based Modelling and System Identification.

\bibitem{Bulea00}
B.~{Bulut}, M.~R. {Katebi}, and M.~A. {Johnson}, ``Industrial application of
  model based predictive control as a supervisory system,'' in {\em Proceedings
  of the 2000 American Control Conference}, vol.~6, pp.~3763--3767 vol.6, June
  2000.

\bibitem{OveDem12}
P.~Van~Overschee and B.~De~Moor, {\em Subspace identification for linear
  systems: Theory—Implementation—Applications}.
\newblock Springer Science \& Business Media, 2012.

\bibitem{Bil95}
P.~Billingsley, {\em Probability and Measure}.
\newblock Wiley Series in Probability and Statistics, Wiley, 1995.

\bibitem{Bemea02}
A.~Bemporad, M.~Morari, V.~Dua, and E.~N. Pistikopoulos, ``The explicit linear
  quadratic regulator for constrained systems,'' {\em Automatica}, vol.~38,
  no.~1, pp.~3 -- 20, 2002.

\end{thebibliography}
\bibliographystyle{ieeetr}

\end{document}